\documentclass[10pt, twoside, a4paper]{article}
\let\finishall\relax\let\Finishall\relax\let\getprepared\relax
\ifx\TestIngCommand\undefined\relax
\else\input ../../../proc/ppreamb-book.tex 
  \input init.tex \fi
\let\TestIngCommand\undefined

\newtheorem{remark}{Remark}

\usepackage{amsmath,amscd,amssymb}                                         
\usepackage{graphics}                                                     
\newtheorem{theo}{Theorem}                                                 
\newtheorem{lem}{Lemma}                                                    
\newtheorem{cor}{Corollary}                                                
\newskip\ttglue\ttglue=.5em plus.25em minus.15em                           
\def\firstname#1{\def\FIRSTNAME{#1}\ignorespaces}
\def\lastname#1{\def\LASTNAME{#1}\ignorespaces}
\def\middleinitial#1{\def\MIDDLEINI{#1}\ignorespaces}
\def\department#1{\def\DEPARTMENT{#1}\ignorespaces}
\def\institute#1{\def\INSTITUTE{#1}\ignorespaces}
\def\address#1{\def\ADDRESS{#1}\ignorespaces}
\def\country#1{\def\COUNTRY{#1}\ignorespaces}
\def\otheraffiliation#1{\def\OTHERAFFILIATION{#1}\ignorespaces}
\def\email#1{\def\EMAIL{#1}\ignorespaces}
\newcount\autcount\autcount=0                                              
\newcount\affcount\affcount=0                                              
\newcount\numcount\newcount\nummcount                                      
\newcount\nummmcount\newcount\nummmmcount                                  
\def\writename#1#2{\ \kern-1ex\hbox{
  \csname AUthor\the#1\endcsname\                                          
  \edef\TESTSTR{}\expandafter\ifx\csname auTHor\the#1\endcsname\TESTSTR    
  \else\csname auTHor\the#1\endcsname.\ \fi                                
  \csname authOR\the#1\endcsname$^{\csname AFF\the#1\endcsname}$
  \expandafter\ifx\csname corr\number#1\endcsname\relax                    
  \else\thanks{Corresponding author.}\ \fi                                 
  }\ifnum#1<#2, \else\ \kern-1ex\fi}
\def\writeemail#1{
  \nummcount=0\relax\nummmcount=0\relax                                    
  \loop\ifnum\nummcount<\autcount\advance\nummcount by1\relax              
    {\expandafter\ifnum\csname AFF\the\nummcount\endcsname=#1\relax        
    \global\advance\nummmcount by1\fi}\repeat                              
  \nummcount=0\relax\nummmmcount=0\relax                                   
  \loop\ifnum\nummcount<\autcount\advance\nummcount by1\relax              
    {\expandafter\ifnum\csname AFF\the\nummcount\endcsname=#1\relax        
    \global\advance\nummmmcount by1\relax\def\blank{}\expandafter          
    \ifx\csname EMAIL\the\nummcount\endcsname\blank(no e-mail)
    \else\csname EMAIL\the\nummcount\endcsname                             
    \fi                                                                    
    \ifnum\nummmmcount<\nummmcount; \fi\fi}\repeat}
\long\def\BeginAuthorList#1\EndAuthorList{#1\relax                         
  \author{\vbox{\hsize=390pt\noindent\numcount=0\relax                     
    \loop\ifnum\numcount<\autcount\advance\numcount by1\relax              
      \writename{\numcount}{\autcount}
      \repeat}\\[2mm]                                                      
    \vbox{\small\numcount=0\relax                                          
      \loop\ifnum\numcount<\affcount\advance\numcount by1\relax            
        \vbox{{\count0=\numcount\relax                                     
          \loop\expandafter\ifnum\csname AFF\the\count0\endcsname
            <\numcount\relax\advance\count0 by1\relax\repeat               
          $^{\csname AFF\the\count0\endcsname}$}
        \def\BLANK{}\expandafter\ifx\csname DEPT\the\numcount\endcsname    
          \BLANK                                                           
          \else\csname DEPT\the\numcount\endcsname, \fi                    
        \csname INST\the\numcount\endcsname,                               
        \csname ADDR\the\numcount\endcsname,                               
        \csname COUN\the\numcount\endcsname                                
        \edef\TEST{}\expandafter\ifx\csname OTHE\the\numcount\endcsname
          \TEST                                                            
          .\else;\break\csname OTHE\the\numcount\endcsname.\fi}
        \vbox{\writeemail{\numcount}}
        \repeat}\\}}
\expandafter\def\csname x1\endcsname{}
\expandafter\def\csname x2\endcsname{}
\expandafter\def\csname x3\endcsname{}
\expandafter\def\csname x4\endcsname{}
\expandafter\def\csname x5\endcsname{}
\expandafter\def\csname x6\endcsname{}
\expandafter\def\csname x7\endcsname{}
\expandafter\def\csname x8\endcsname{}
\expandafter\def\csname x9\endcsname{}
\def\Author#1#2{\global\advance\autcount by1\relax#2                       
  \expandafter\edef\csname AUthor\the\autcount\endcsname{\FIRSTNAME}
  \expandafter\edef\csname auTHor\the\autcount\endcsname{\MIDDLEINI}
  \expandafter\edef\csname authOR\the\autcount\endcsname{\LASTNAME}
  \expandafter\edef\csname EMAIL\the\autcount\endcsname{\EMAIL}
  \let\tempera\"\def\"{\string\"}\expandafter\ifx\csname x\DEPARTMENT
    \endcsname\relax                                                       
    \global\advance\affcount by1\relax\let\"\tempera                       
    \expandafter\edef\csname DEPT\the\affcount\endcsname{\DEPARTMENT}
    \expandafter\edef\csname INST\the\affcount\endcsname{\INSTITUTE}
    \expandafter\edef\csname ADDR\the\affcount\endcsname{\ADDRESS}
    \expandafter\edef\csname COUN\the\affcount\endcsname{\COUNTRY}
    \expandafter\edef\csname OTHE\the\affcount\endcsname{\OTHERAFFILIATION}
    \expandafter\edef\csname AFF\the\autcount\endcsname{\the\affcount}
  \else\expandafter\edef\csname AFF\the\autcount\endcsname{\DEPARTMENT}
  \fi\let\"\tempera\ignorespaces}
\def\CorrespondingAuthor#1#2{
  \expandafter\xdef\csname corr\number#1\endcsname{cor}
  \Author#1{#2}}
\def\PaperTitle#1{\title{\bf#1}}
\def\Category#1{\ignorespaces}
\def\keywords#1{{\noindent \emph{Keywords:}                                
  \def\BLANK{}\def\TEST{#1}\ifx\BLANK\TEST(n/a).\else#1\fi}}
\getprepared                                                               
\begin{document}                                                           

\PaperTitle{Spectral Factorization and Entire Functions}
\Category{(Pure) Mathematics}

\date{}

\BeginAuthorList
  \Author1{
    \firstname{Wayne}
    \lastname{Lawton}
    \middleinitial{M}   
    \department{Department of the Theory of Functions, Institute of Mathematics and Computer Science}
    \institute{Siberian Federal University}
    \otheraffiliation{}
    \address{Krasnoyarsk}
    \country{Russian Federation}
    \email{wlawton50@gmail.com}}
\EndAuthorList
\maketitle
\thispagestyle{empty}
\begin{abstract}
The Fej\'{e}r-Riesz spectral factorization lemma, which represents a nonnegative trigonometric polynomial as the squared modulus of a trigonometric polynomial, was extended by Ahiezer to factorize certain entire functions and by Helson and Lowdenslager to factorize certain functions on compact connected abelian groups whose Pontryagin duals are equipped with a linear order. This paper relates these factorizations for archimedian orders using the theory of almost periodic functions.
\end{abstract}

\finishall
\section{Introduction}
$\mathbb N = \{1,2,3,...\},$ $\mathbb Z,$ $\mathbb Q,$ $\mathbb R,$ $\mathbb C$
are the natural, integer, rational, real, and complex numbers.
$\mathbb D := \{ \, z \in \mathbb{C} \, : \, |z| \leq 1 \, \}$
is the closed unit disk, $\mathbb D^o$ is its interior and its boundary
$\mathbb T$ is the multiplicative circle group. For $z \in \mathbb C,$ define
$\chi_z \, : \, \mathbb R \rightarrow \mathbb C$ by $\chi_z(s) = e^{izs}.$ \\ \\
In 1915 Fej\'{e}r conjectured \cite{fejer} and Riesz proved \cite{riesz1} that if $f$ is a nonnegative trigonometric polynomial of the form
$f(t) = \sum_{-d}^{d} f_k \, e^{2\pi i k t}$
where $f_d \neq 0,$ then there exists a trigonometric polynomial
$h(t) = \sum_{0}^{d} h_k \, e^{i 2\pi k t}$
such that
$f = |h|^2.$
Moreover, there exists a unique such $h$ so that the polynomial
    $H(z) := \sum_{0}^{d} h_k \, z^k$
has no zeros in $\mathbb D^o$ and $H(0) > 0.$ This result is called the Fej\'{e}r-Riesz spectral factorization lemma. Riesz gave a proof based on the fundamental theorem of algebra in (\cite{riesz2}, p. 117).
\\ \\
In 1921 Szeg\"{o} \cite{szego1} extended the Fej\'{e}r-Riesz lemma by proving that if $w \in L^1(\mathbb T)$ is nonnegative then $\log w \in L^1(\mathbb T)$ iff there exists a nonzero function $h$ in the Hardy space $H^2(\mathbb T)$ such that $w = |h|^2.$ Moreover $h$ is unique if its
holomorphic extension $H$ to $\mathbb D^o$ has no zeros and $H(0) > 0.$
$H^p(\mathbb T), p \in [1,\infty]$ consist of functions in $L^p(\mathbb T)$ that are nontangential boundary values of their holomorphic extensions to $\mathbb D^o$ (\cite{rudin2}, Theorem 17.11).
\\ \\
Let ${\mathcal{E}}(\mathbb{R})$ be the Fr\'{e}chet space of
smooth functions with the topology of uniform convergence of derivatives
over compact subsets, ${\mathcal{E}}^{\prime}(\mathbb{R})$ be its
dual space of compactly supported distributions \cite{gelfand2, schwartz1},
and $<\cdot,\cdot> \, : \, {\mathcal{E}}^{\prime}(\mathbb{R}) \times {\mathcal{E}}(\mathbb{R}) \rightarrow \mathbb C$
be the bilinear pairing. The support interval of $\eta \in {\mathcal{E}}^{\prime}(\mathbb{R})$ 
is the smallest closed interval $[\alpha,\beta]$ containing its support. Then
$F(z) := \, <\eta,\chi_z>$
is entire of exponential type
$\tau := \limsup_{|z| \rightarrow \infty} |z|^{-1} \log |F(z)| = \max \{|\alpha|,|\beta|\}$
and its restriction to $\mathbb R$ is bounded by a polynomial.
The Paley-Wiener-Schwartz theorem \cite{paley}, \cite{schwartz2}
implies the converse.
\\ \\
In 1948 Ahiezer \cite{ahiezer} proved that if $F \, : \, \mathbb C \rightarrow \mathbb C$ is an entire function of exponential type $\tau > 0,$ the restriction $F|_{\mathbb R} \geq 0,$ and
$\int_{-\infty}^\infty (1+x^2)^{-1}\, \max\{\, 0, \log F(x)\, \}\, dx < \infty,$
then there exists an entire function $S$ of exponential type $\tau/2$ with no zeros in $\mathbb U^o$  such that
$F(z) = S(z)\, \overline {S(\overline z)}$
and $S$ is unique up to multiplication by a constant with modulus $1.$
Boas (\cite{boas2}, Theorem 7.5.1) gives a proof based on the Ahlfors-Heins theorem \cite{heins} and proves that Ahiezer's theorem implies the Fej\'{e}r-Riesz lemma. \\ \\
Henceforth assume $G$ is a compact abelian topological group with a linear order $\leq$ on its Pontryagin dual \cite{pontryagin} $\widehat G,$ which consists of characters or continuous homomorphisms $\chi : G \rightarrow \mathbb T,$ under pointwise multiplication.
Such an order exists iff the identity $1 \in \widehat G$ is the only element of finite order (\cite{rudin1}, Theorem 8.1.2) iff $G$ is connected \cite{levi}.
The Stone-Weierstrass theorem (\cite{rudin2}, A14) implies that the algebra of trigonometric polynomials $T(G) := \hbox{span } \widehat G$ is dense in the $C^*$-algebra $C(G).$
Let $d\sigma$ be Haar measure on $G$ normalized so $\int_G d\sigma = 1,$
and let $L^p(G), \, p \in [1,\infty]$ be the associated Banach spaces and
$L_r^p(G) := \{ \, f \in L^p(G) \, : \, f(G) \subset \mathbb R \, \}.$
For $p < \infty$ and $f \in L^p(G)$ there exists a sequence of nonnegative $e_n \in T(G),$ with $\int_G e_n d\sigma = 1,$
such that the sequence of convolutions
$e_n * f \in T(G)$
satisfies
$\lim_{n \rightarrow \infty} ||f-e_n * f||_p = 0$
hence by (\cite{rudin2}, Theorem 3.12) there exists a subsequence of $e_n*f$ that converges pointwise to $f.$ Let $L_{\ast}^1(G)$ be the subset of $L^1(G)$ consisting of functions for which there exist such a sequence $e_n$ with
\begin{equation}\label{sup}
    \sup_n |e_n*f| \in L^1(G).
\end{equation}
We record that $L_{\ast}^1(\mathbb T)$ is the subset of $f \in L^1(\mathbb T)$
whose Hardy-Little maximal function
\begin{equation}\label{max}
    (Mf)(z) := \sup_{L \in [0,\pi]} \frac{1}{2L} \int_{-L}^L f(ze^{is}) \, ds, \ \ z \in \mathbb T
\end{equation}
is in $L^1(\mathbb T)$ (\cite{laugesen}, Theorem 7.3). The Fourier transform $\mathfrak{F} \, : \, L^1(G) \rightarrow \ell^\infty(\widehat G)$ is defined by
$\mathfrak{F}(h)(\chi) \mathrel{\mathop :}= \int_G \, h \, \overline \chi \, d\sigma.$
The Hausdorff-Young inequality \cite{hausdorff, young} implies that for $p \in [1,2]$ and $p^{-1} + q^{-1} = 1$ its restriction is a bounded operator
$\mathfrak{F} \, : \, L^p(G) \rightarrow \ell^q(\widehat G).$
The spectrum of $h \in L^1(G)$ is $\Omega(h) := \hbox{ support } \mathfrak{F}(h).$
Define $HT(G) = \hbox{span } \{\, \chi \in \widehat G \, : \, 1 \leq \chi \, \}$
and the Hardy space
$H^p(G) := \{\, h \in L^p(G) \, : \, \chi < 1 \Rightarrow \chi \notin \Omega(h) \, \}.$
\\ \\
The Mahler measure \cite{mahler1}, \cite{mahler2} of $h \in L^1(G)$ is
\begin{equation}
\label{Mahler}
M(h) := \exp \int_{G} \log |h| \, d\sigma.
\end{equation}
Clearly $M(h) > 0$ iff $\log |h| \in L^1(G)$ and $M(h) = 0$ otherwise. $h \in L^1(G)$ is outer if $M(h) = |\int_G h \, |$ and inner if $|h| = 1$ almost everywhere.
\\ \\
In 1949 Beurling \cite{beurling} introduced these concepts for $G = \mathbb T$ and proved that every nonzero $h \in H^p(\mathbb T)$ satisfies $\log |h| \in L^1(\mathbb T)$ and admits a unique factorization $h = h_i\, h_o$ where $h_i \in H^{\infty}(\mathbb T)$ is inner and $h_o \in H^p(\mathbb T)$ is outer and $\int_{\mathbb T} h_o > 0.$ The restriction $h$ of an algebraic polynomial $H$ to $\mathbb T$ is in $H^1(\mathbb T)$ and is outer iff $H$ has no roots in $\mathbb D^o$ since Jensen's formula \cite{jensen} gives $M(h) = |h(0)| \prod_{h(\lambda) = 0} \max \{|\lambda|^{-1},1\}.$ Jensen's formula also implies that $h \in H^p(\mathbb T)$ is outer iff its holomorphic extension $H$ to $\mathbb D^o$ has no zeros.
\\ \\
In 1958 Helson and Lowdenslager extended results of Szeg\"{o} and Beurling to groups with ordered duals. Rudin in (\cite{rudin1}, Theorem 8.4.3) explains their main result: if $w \in L^1(G)$ is nonnegative then $\log w \in L^1(G)$ iff there exists a $h \in H^2(G)$ with $\int h \neq 0$ such that $w = |h|^2.$ Then there is a unique such outer $h$ satisfying $\int h > 0$ which we denote
by $S(w)$ and call the spectral factor of $w.$ This result implies an extension of Beurlings factorization result. If $h \in H^2(G)$ then $w = |h|^2 \in L^1(G)$ and Rudin (\cite{rudin2}, Theorem 8.4.1) showed that if $\int_G h \neq 0$ then $\log w \in L^1(G)$ and $h = h_ih_o,$ where $h_o$ is the spectral factor of $w$ and
$h_i := h_o^{-1}h,$ is the Beurling factorization of $h.$. \\ \\
In Section \ref{almost} we use Ahiezer's spectral factorization to study Hardy spaces of almost periodic functions. 
In Section \ref{spectral} we prove, under slightly more restrictive conditions on $w,$ that its Helson-Lowdenslager spectral factor can be expressed as $h = \exp v$ where $v \in H^2(G).$ In Section
\ref{compact} we relate the two factorizations when the order on $\widehat G$ is archimedian.
Then the order gives a homomorphism $\theta : \mathbb R \rightarrow G$ having a dense image
and the map $\Theta \, : \, C(G) \rightarrow U(\mathbb R),$ defined by $\Theta(f) := f \circ \theta,$ enables us to related these factorizations. Our main result is that if $\Omega(w)$ is bounded then $\Omega(h)$ is bounded.
\section{Ahiezer's Factorization and Almost Periodicicity}\label{almost}
$\mathbb U := \{\, z \in \mathbb C \, : \, \Im z \geq 0 \, \}$ and $\mathbb U^o :=$ its interior.
$HAR(\mathbb U); HOL(\mathbb U)$ is the set of continuous $\mathbb C$-valued functions on $\mathbb U$ whose restrictions to $\mathbb U^o$ are harmonic; holomorphic, respectively. $HOL(\mathbb U) \subset HAR(\mathbb R)$ and a harmonic function $F$ is holomorphic iff $\Re F$ and $\Im F$ satisfy the Cauchy-Riemann equations.
$EE_{\tau}(\mathbb C) :=$ set of entire functions of exponential type $\tau \geq 0.$ \\ \\
Characters on $\mathbb R$ have the form
$\chi_\omega$ with $\omega \in \mathbb R.$ $T(\mathbb R)$ is the algebra of trigonometric polynomials on $\mathbb R$ spanned by characters. $C_b(\mathbb R)$ is the $\textrm C^*$-algebra of bounded continuous functions with norm $||f|| = \sup_{x \in \mathbb R} |f(x)|.$ Bohr \cite{bohr1} defined the $\textrm C^*$-algebra $U(\mathbb R)$ of uniformly almost periodic functions to be the closure of $\textrm{T}(\mathbb R)$ in $\textrm C_b(\mathbb R)$ and proved that their means
$m(f) \, \mathrel{\mathop :}= \lim_{L \rightarrow \infty} (2L)^{-1} \int_{-L}^L f(t) dt$
exist.
The Fourier transform of $f \in U(\mathbb R)$ is
$\widehat f(\chi) := m(f\, \overline \chi),$ and the spectrum is
$\Omega(f) := \hbox{ support } \widehat(f).$ Parseval's identity
$m(|f|^2) = \sum_{\omega \in \Omega(f)} |\widehat f(\omega)|^2$
holds hence $\Omega(f)$ is countable.
A nonzero function $f \in U(\mathbb R)$ has a bounded spectrum if its bandwidth
$b(f) := \sup \Omega(f) - \inf \Omega(f) < \infty.$ We define
$$AU(\mathbb R) := \{\,f \in U(\mathbb R) \, : \, {\widehat f} \in \ell^1(\Omega(f)) \, \},\ \
BU(\mathbb R) := \{\,f \in U(\mathbb R) \, : \, b(f) < \infty \, \},$$
$$HU(\mathbb R) := \{\, h \in HU(\mathbb R) \, : \, \Omega(f) \subset [0,\infty) \, \}, \ \
IHU(\mathbb R) := \{\, h \in HU(\mathbb R) \, : \, h^{-1} \in HU(\mathbb R) \, \}.$$
%
%
\begin{lem}\label{lem1}
If $f \in U(\mathbb R)$ then $f \in BU(\mathbb R)$ iff
$f = F_{\mathbb R}$ for some $F \in EE_\tau(\mathbb C)$ where $\tau = \max \{\, |\inf \Omega(f)|, \, |\sup \Omega(f)|\, \}.$
\end{lem}
{\bf Proof} Bohr \cite{bohr2} first proved this. Also see Szeg\"{o} \cite{szego2}.
English language proofs were given by Bohr (\cite{bohr3}, Appendix II, p. 114--115) and by Levin (\cite{levin}, p. 268, Corollary).
\begin{remark}\label{rem1}
If $f \in AU(\mathbb R) \cap BU(\mathbb R)$ then
$F(z) = \sum_{\omega \in \Omega(f)} \widehat f(\omega)\, e^{i\omega z}$
since the sum converges absolutely uniformly over compact subsets.
\end{remark}
The Poisson kernel functions $P_y \, : \, \mathbb R \rightarrow \mathbb R, y > 0$ are
$P_y(x) := \frac{y}{\pi} \frac{1}{x^2 + y^2}, \ x + iy \in \mathbb U^o.$
The Poisson integral $H \, : \, \mathbb U \rightarrow \mathbb C$
of $h \in C_b(\mathbb R)$ is defined by $H(x) := h(x), x \in \mathbb R$ and
\begin{equation}\label{Poissonint}
    F(x+iy) := (P_y * f)(x) = \frac{y}{\pi} \int_{-\infty}^{\infty} \frac{f(x-s)}{s^2 + y^2} \,  ds,  \ \ x+iy \in \mathbb U^o.
\end{equation}
%
%
\begin{lem}\label{lem2}
If $h \in U(\mathbb R)$ then its Poisson integral $H \in HAR(\mathbb U)$
satisfies \newline
$\lim_{y \rightarrow \infty} \sup_{x \in \mathbb R} |H(x+iy) - m(h)| = 0,$ and
$h \in HU(\mathbb R)$ iff $H \in HOL(\mathbb U).$
\end{lem}
{\bf Proof}
The first assertion is standard, the second since
$\{P_y:y > 0\}$ is an approximate identity and $aP_y(ax) = P_y(a^{-1}y), \, a > 0,$
the third since the Poisson integral of $\chi_\omega$ equals $e^{i\omega\, z}$ for $\omega \geq 0$ and
equals $e^{i\omega\, \overline z}$ for $\omega < 0.$
%
%
\begin{remark}\label{rem2}
If $f \in AU(\mathbb R)$ then
$F(z) = \sum_{\omega < 0} \widehat f(\omega)\, e^{i\omega \overline z} + \widehat f(0) +
\sum_{0 < \omega} \widehat f(\omega)\, e^{i\omega z}$ since the sum converges absolutely uniformly over compact subsets of $\mathbb U.$
\end{remark}
%
%
\begin{lem}\label{lem3}
If $h \in IHU(\mathbb R)$ then $\widehat h(0) \neq 0,$ the Poisson integral $H$ of $h$ has no zeros, and $|H|$ is bounded below by a positive number. Therefore $\Re \log H = \log |H|$ is bounded.
\end{lem}
{\bf Proof}
Since $hh^{-1} = 1$ implies $1 = \sum_{\omega_1 + \omega_1 = 0} \widehat h(\omega_1) \widehat {h^{-1}}(\omega_2) = \widehat h(0) \widehat {h^{-1}}(0)$ hence $\widehat h(0) \neq 0.$
Let $H$ and $H_1$ be the Poisson integrals of $h$ and $h^{-1},$ respectively
and define $J(z) := H(z)H_1(z), \Im z \geq 0$ and
$J(z) := \overline {H(\overline z)H_1(\overline z)}, \Im z \leq 0.$
The Schwarz reflection principle (\cite{rudin2}, Theorem 11.14 ) implies that $J$ is entire. Since it is bounded Liouville's theorem (\cite{rudin2}, Theorem 10.23) implies that $J = 1$ hence $H$ has no zeros. Since $|H_1|$ is bounded above, $|H|$ is bounded below by $\left[\, \sup_{z \in \mathbb U} |H_1(z)|\, \right]^{-1} > 0.$ Since $H$ has no zeros, $\log \, H$ exists and the conclusion follows directly.
%
%
\begin{lem}\label{lem4}
$E \in HOL(\mathbb U), \Re \, E$ is bounded above, and $\Re \, E|_{\mathbb R} = 0$ then $E(z) = \alpha i + \beta iz$ where $\alpha \in \mathbb R$ and $\beta \geq 0.$
\end{lem}
{\bf Proof} Deng \cite{deng} proved that if $U \in HAR(\mathbb U)$ is real valued and
$U^{+} := \max \{\, U, 0 \, \}$ satisfies mild growth conditions then $U$ is a boundary integral of $u := U|_{\mathbb R}$ plus a harmonic polynomial that vanishes on $\mathbb R.$ Therefore $E$ is a polynomial having the form above.
%
%
\begin{theo}\label{thm1}
If $h \in IHU(\mathbb R)$ and $f := |h|^2 \in BU(\mathbb R)$ then
$h(x) = e^{-\alpha \, i}\, e^{i\frac{\tau}{2}x}\, S(x)$ where $\alpha \in \mathbb R,$
$\tau = b(f)/2$ and $S$ is Ahiezer's spectral factor of the entire function $F$ such that $F|_{\mathbb R} = f.$
Therefore $h, S \in BU(\mathbb R),$ $\Omega(h) \subset [0,\tau]$ and
$\Omega(S) \subset [-\frac{\tau}{2},\frac{\tau}{2}].$
\end{theo}
{\bf Proof}
Lemma \ref{lem1} implies there exists $F \in EE_{\tau}(\mathbb R)$ with $\tau = b(f)/2$ and
$F|_{\mathbb R} = f.$
Then $F$ satisfies the hypothesis of Ahiezer's theorem \cite{ahiezer} so there exists $S \in EE_{\tau/2}(\mathbb C)$
such that $S|_{\mathbb U^o}$ has no zeros and $F(z) = S(z)\overline {S(\overline z)}.$ Define $S_1(z) := e^{i\frac{\tau}{2}z}\, S(z)|_{\mathbb U}.$ Then $S_1 \in HOL(\mathbb U)$ has no zeros and $|S_1|$ is bounded above. Therefore $E_1 := \log S_1$ exists and $\Re\, E_1 = \log |S_1|$ is bounded above.
Let $H \in HOL(\mathbb U)$ be the Poisson integral of $h.$
Lemma \ref{lem3} ensures that $E := \log H$ exists and that $\Re \, E = \log |H|$ is bounded above and below.
Therefore $E_1 - E \in HOL(\mathbb U),$ $\Re (E_1 - E)$ is bounded above, and
$\Re (E_1-E)|_{\mathbb R} = 0$ so Lemma \ref{lem3} implies $E_1(z) = E(z) + \alpha i + \beta i z, \alpha \in \mathbb R, \beta \geq 0$ so
$S_1(z) = e^{\alpha i + \beta i z}H(z).$ If $\beta > 0$ then
$S$ would have exponential type $< \frac{\tau}{2}$ thus contradicting Ahiezer's theorem.
Therefore $h(x) = e^{-\alpha \, i}S_1(x) = e^{-\alpha \, i}\, e^{i\frac{\tau}{2}x}\, S(x).$
\section{Helson-Lowdenslager Spectral Factorization}\label{spectral}
If $f$ is a function on $G$ and $f_n$ is a sequence of functions on $G$ we let
$f_n \rightarrow f$ denote pointwise convergence almost everywhere with respect to Haar measure.
Define $\widehat G_+ := \{\, \chi \in \widehat G \, : \, \chi > 1 \, \}$ and
$T_+(G) := \{\, P \in \widehat T(G) \, : \, \Omega(P) \subset \widehat G_+ \, \}.$
Assume that $w \in L^1(G),$ $w \neq 0,$ and $w \geq 0.$ Let
$L_w^2(G)$ be the Hilbert space completion of $C(G)$ with the norm
$||f||_w := \int_G |f|^2 \, w \, d\sigma,$
and let $K$ be the closure of $1 + T_+(G)$ in $L_w^2(G).$
%
%
\begin{lem}\label{lem5}
There exists a unique $H$ in the closure of $T_+(G)$ such that
\begin{equation}
    \int_G |1 + H|^2 \, w \, d\sigma = \inf_{P \in T_+(G)} \int_G |1 + P|^2 \, w \, d\sigma.
\end{equation}
$|1+H|^2 w$ is constant a.e. with respect to $d\sigma.$
Define $h := \sqrt {M(w)} \, (1 + H)^{-1}.$
Then
\begin{equation}\label{hinth}
     h \in H^2(G) \ \ \hbox{and} \ \ \int_G h = \sqrt {M(w)}.
\end{equation}
Therefore $w = |h|^2$ and $h$ is outer so $h = S(w).$
\end{lem}
{\bf Proof}
The first statement follows since $K$ is nonempty, closed and convex so it contains a unique element $1+H$ having minimum norm (\cite{rudin2}, Theorem 4.10), then the second follows since
$\chi \in \widehat G_+ \Rightarrow K(1+\chi) \subset K$
and hence
$\int_G |1 + H|^2 \, w \, \chi \, d\sigma = 0$ for all $1 \neq \chi \in \widehat G.$
The third assertion is deeper and the fourth is obvious. All four assertions were proved by Szeg\"{o} \cite{szego1} for $G = \mathbb T$ and by Helson and Lowsenslager in
(\cite{helsonlow1}, Theorem 1, Lemma 2, Theorem 3) for scalar value functions on $\mathbb T^2$ and in (\cite{helsonlow1}, Section 6) for matrix valued functions on general groups. Rudin
(\cite{rudin1}, Chapter 8) derives the third assertion for scalar valued functions on general groups.
\begin{lem}\label{lem6}
If $w \in L^1(\mathbb T)$ is nonnegative and $\log w \in L^1(\mathbb T)$ then
$$S(w)(z) = \lim_{n \rightarrow \infty} \exp \left[\frac{1}{2\pi} \int_{-\pi}^{\pi}
\frac{e^{it}+ (1-n^{-1})\, z}{e^{it} - (1-n^{-1})\, z}\, \frac{\log w}{2}\, dt \right].$$
\end{lem}
{\bf Proof} See (\cite{rudin2}, Theorem 17.16). \\ \\
Lemma \ref{lem6} expresses $S(w)$ as a limit of a sequence of invertible functions in $H^2(\mathbb T).$
We now derive a similar result for functions on more general groups. The Hilbert transform $\mathcal{H} :  L_r^2(G) \rightarrow L_r^2(G)$ and analytic transform
$\mathcal{A} :  L_r^2(G) \rightarrow H^2(G)$
\begin{equation}\label{transforms}
    \mathfrak{F}(\mathcal{H}(f))(\chi) := -i \, sign(\chi)\,  \mathfrak{F}(f)(\chi) \, ; \ \ \ \
    \mathfrak{F}(\mathcal{A}(f)):= f + i \mathcal{H}(f)
\end{equation}
are bounded, $\mathcal{A}(1) = 1,$ $\mathcal{A}(\chi) = 0$ if $\chi < 1,$
$\mathcal{A}(\chi) = 2\, \chi$ if $\chi > 1,$ and $\Re \mathcal{A}(f) = f, f \in L_r^2(G).$
The Wiener algebra $A(G) := \{ \, f \in L^1(G)\, : \, \mathfrak{f} \in \ell^1(\widehat G) \, \}$ is a Banach algebra under pointwise multiplication and norm $||f||_{A(G)} := ||\widehat f||_{\ell^1(\widehat G)}.$ (\cite{rudin2}, Theorem 5.12) implies that $A(G)$ is a proper subset of $C(G).$ We define
    $A_{\neq 0}(G) := \{ f \in A(G) :  0 \notin f(G) \},$
    $A_r(G) := \{ f \in A(G) : f(G) \subset \mathbb R \},$
    $A_+(G) := \{ f \in A(G) :  f > 0 \},$
    $HA := A(G) \cap H^\infty(G),$
	$OHA(G) := \{ f \in HA(G) : f \hbox{ is outer} \}$ and
    $IOHA(G) := \{ f \in OHA(G) : f^{-1} \in HA(G) \}.$
%
%
%
\begin{lem}\label{lem7}
\begin{enumerate}
\item $f \in A_{\neq 0}(G) \Rightarrow f^{-1} \in A_{\neq 0}(G).$
\item $f \in A(G) \Rightarrow \exp f \in A_{\neq 0}(G)$ and $||\exp f||_{A(G)} \leq \exp ||f||_{A(G)}.$
\item $f \in A_+(G) \Rightarrow \log f \in A_r(G).$
\item $\mathcal{A} \, : \, A_r(G) \rightarrow HA(G)$ is a bounded operator in the $A(G)$-norm.
\item If $\int_G f \, d\sigma \in \mathbb R$ and $f \in HA(G)$
 then $\exp f \in IOHA(G).$
\end{enumerate}
\end{lem}
{\bf Proof}
1. follows from the Gelfand representation \cite{gelfand1} (for $G = \mathbb{T}$ it follows from Weiner's Tauberian lemma \cite{wiener}), 2. and 3. follow from the Arens-Royden theorem \cite{arens}, \cite{royden} and 4. is obvious. Let $h := \exp f.$ Clearly $h, h^{-1} \in HA(G)$ and $h$ is outer since
$$\int_G \exp f \, d\sigma = \exp \int_G f \, d\sigma =
\exp \int_G \Re f \, d\sigma = \exp \int_G \log |h| \,  d\sigma = M(h).$$
Henceforth $w : G \rightarrow [0,\infty)$ is nonzero and measurable,
$u := \frac{1}{2}\log w, \ v := \mathcal{A}(u), \ h = \exp v.$
%
%
\begin{lem}\label{lem8}
$w \in A_+(G) \Rightarrow S(w) = h \in IOHA(G).$
\end{lem}
{\bf Proof}
Lemma \ref{lem7} gives
$w \in A_+(G) \Rightarrow u \in A_r(G) \Rightarrow v \in HA(G) \Rightarrow h \in IOHA(G)$
since $\int_G v \, d\sigma = \widehat v(1) = \widehat u(1) = \int_G u \, d\sigma \in \mathbb R.$
%
%
\begin{lem}\label{lem9}
If $0 < a \leq w \leq b$ then $S(w) = h$ and there exists a sequence
$h_n \in IOHA(G)$ such that $\lim_{n \rightarrow \infty} ||h-h_n||_2 = 0$
and $\Omega(|h_n|^2) \subset \Omega(w).$
\end{lem}
{\bf Proof}
Choose a sequence $e_n \in T(G)$ so that $w_n : = e_n*w \rightarrow w.$
Then $w_n(G) \subset [a,b]$ hence  $w_n \in A_+(G).$
Define $u_n := \frac{1}{2}\log w_n,$ $v_n := \mathcal{A}(u_n),$ and $h_n := \exp v_n.$
Lemma \ref{lem8} implies that $h_n = S(w_n)$ and $h_n \in IOHA_{\neq 0}(G).$
Since $\sup_{n} |u-u_n|^2 \leq \infty$
and $|u-u_n|^2 \rightarrow 0,$
Lebesque's DCT (Dominated Convergence Theorem) (\cite{rudin2}, Theorem 1.34) implies that
$||u-u_n||_2 \rightarrow 0.$
Since $\mathcal{A} : L_r^2(G) \rightarrow H^2(G)$ is bounded,
$\lim_{n \rightarrow \infty} ||v-v_n||_2 = 0,$
so (\cite{rudin2},Theorem 3.12) implies that we could have chosen $e_n$ so that
$|v-v_n| \rightarrow 0$ hence
$|h - h_n| \rightarrow 0.$
Since $|h-h_n|^2 \leq 4b,$ Lebesque's DCT implies that
$\lim_{n \rightarrow \infty} ||h-h_n||_2 = 0.$
Each $h_n \in HA_{\neq 0}(G)$ is outer so $h \in  H^2(G)$ is outer. Clearly
$\Omega(|h_n|^2) = \Omega(e_n*w) \subset \Omega(w).$
%
%
\begin{lem}\label{lem10}
If $w \in L_*^1(G)$ and $0 < a \leq w$ then $S(w) = h$ and there exists a sequence
$h_n \in IOHA(G)$ such that $\lim_{n \rightarrow \infty} ||h-h_n||_2 = 0$
and $\Omega(|h_n|^2) \subset \Omega(w).$
\end{lem}
{\bf Proof}
$u \in L_r^2(G),$ $v \in H^2(G),$ and $|h|^2 = w.$
Choose a sequence $e_n \in T(G)$ with $w_n : = e_n*w \rightarrow w$
and $\sup_{n} |w_n| \in L^1(G).$ Since $a \leq w_n,$ $w_n \in A_+(G).$
Define $u_n := \frac{1}{2}\log w_n,$ $v_n := \mathcal{A}(u_n),$ and $h_n := \exp v_n.$
Lemma \ref{lem8} gives $S(w_n) = h_n \in IOHA(G)$ and $\Omega(|h_n|^2) \subset \Omega(w).$
Since $|u - u_n|^2 \rightarrow 0$
and $\sup_{n} |u - u_n|^2 \in L^1(G),$
Lebesgue's DCT gives
$\lim_{n \rightarrow \infty} ||u - u_n||_2 \rightarrow 0.$
$\mathcal{A}$ is bounded and
$\lim_{n \rightarrow \infty} ||v - v_n||_2 \rightarrow 0,$ therefore
(\cite{rudin2}, Theorem 3.12) implies that we could have chosen $e_n$ such that
$v_n \rightarrow v.$
Then $|h - h_n|^2 \rightarrow 0.$
Since $\sup_{n} |h - h_n|^2 \in L^1(G),$ Lebesgue's DCT gives
$\lim_{n \rightarrow \infty} ||h - h_n||_2 \rightarrow 0,$
therefore $h = S(w).$
%
%
\begin{theo}\label{thm2}
If $w \in L_*^1(G),$ $w \geq 0,$ and $\log w \in L^2(G)$ then $S(w) = h$ and there exists a sequence
$h_n \in IOHA(G)$ such that $\lim_{n \rightarrow \infty} ||h-h_n||_2 = 0$
and $\Omega(|h_n|^2) \subset \Omega(w).$
\end{theo}
{\bf Proof}
Define $w_m := w + \frac{1}{m}$
$u_m := \frac{1}{2} \log w_n,$ $v_n := \mathcal{A}(u_n), \, f_n := \exp v_n.$
Then $|u - u_n|^2 \rightarrow 0$ monotonically so
Lebesgue's MCT (Monotone Convergence Theorem) (\cite{rudin2}, Theorem 1.26) implies that
$\lim_{n \rightarrow \infty} ||u - u_n||_2 = 0.$
Since $\mathcal{A}$ is bounded and
$\lim_{n \rightarrow \infty} ||v - v_n||_2 = 0,$
(\cite{rudin2},Theorem 3.12) implies that there exist a subsequence $k_n$ of integers with $v_{k_n} \rightarrow v.$ By replacing $w_n$ by $w_{k_n} = w + \frac{1}{k_n}$ we
have $v_n \rightarrow v$ therefore
$|h - f_n|^2 \rightarrow 0$ so Lebesgue's DCT implies
that $\lim_{n \rightarrow \infty} ||h - f_n||_2 = 0.$
Lemma \ref{lem10} implies that for every $n \geq 1,$
$S(w_n) = f_n$ and there exists a sequence $h_{n,m} \in IOHA(G)$ such
that $\lim_{m \rightarrow \infty} ||f_n-h_{n,m}||_2 = 0$
and $\Omega(|h_{n,m}|^2) \subset \Omega(w_n) = \Omega(w)$ since $1 \in \Omega(w).$
Then there exists a sequence $m_n$ such that $h_n := h_{n,m_n}$
satisfies $\lim_{n \rightarrow \infty} ||h - h_n||_2 = 0.$ Then $h \in H^2(G)$ is outer
and $S(w) = h.$
\begin{remark}\label{rem3} $1+z \in HA(\mathbb T)$
is outer but $(1+z) \notin IOHA(\mathbb T).$ However $1+z$ is the limit in $H^2(\mathbb T)$ of the sequence $1 + (1-\frac{1}{n})z \in IOHA(\mathbb T).$
\end{remark}
\section{Compactifications and Archimedean Orders}\label{compact}
A compactification $(G,\theta)$ of $\mathbb R$ consists of a compact connected abelian group $G$ and a continuous homomorphism $\theta \, : \, \mathbb R \rightarrow G$ with a dense image.
We define $\widehat \theta\, : \, \widehat G \rightarrow \mathbb R$ by
\begin{equation}\label{hattheta}
    \widehat \theta \, (\chi) := \omega \hbox{ where } \chi \circ \theta = \chi_\omega
\end{equation}
and an associated order $\leq$ on $\widehat G$ by
\begin{equation}\label{order}
    \chi_1 \leq \chi_2 \hbox{ iff } \widehat \theta(\chi_1) \leq \widehat \theta(\chi_2).
\end{equation}
Since $\theta(R)$ is dense in $G,$ $\widehat \theta$ is injective hence $\leq$ is an archimedian order on $\widehat G.$
\begin{lem}\label{lem11}
Every archimedian order on $\widehat G$ arises from a compactification $(G,\theta)$ as described by Equation \ref{order}. Compactifications $(G,\theta_1)$ and $(G,\theta_2)$ give the same order iff there exists $a > 0$ such that
\begin{equation}\label{thetaequiv}
\theta_2(t) = \theta_1(at), t \in \mathbb R.
\end{equation}
\end{lem}
{\bf Proof}
Otto H\"{o}lder proved  (\cite{rudin1}, Theorem 8.1.2, p. 194), (\cite{ribenboim}, p. 60) that every archimedian order $\leq$ on $\widehat G$ is induced by an injective homomorphism
$\varphi \, : \, \widehat G \rightarrow \mathbb R.$
Define
$\widehat \varphi : \widehat {\mathbb R} \rightarrow \widehat {\widehat G}$
by
$\widehat \varphi(\chi_{\omega}) := \chi_{\omega} \circ \varphi$
and define
$\theta : \mathbb R \rightarrow \widehat {\widehat G}$ by
$\theta(\omega) := \widehat \varphi(\chi_\omega).$
Since $\varphi$ is injective $\widehat \varphi$ and therefore $\theta$ has a dense image.
Pontryagin's duality theorem (\cite{rudin1}, Theorem 1,7.2) implies that $\widehat {\widehat G} = G$ (canonically isomorphic) which proves the first assertion. The second assertion follows since injective homomorphisms $\varphi_1, \varphi_2 : \widehat G \rightarrow \mathbb R$ give the same archimedian order iff there exists $a > 0$ such that $\varphi_2(\chi) = a\, \varphi_1(\chi), \chi \in \widehat G$ hence the corresponding functions $\theta_1, \theta_2$ satisfy
Equation \ref{thetaequiv}. \\ \\
Henceforth we equip $\widehat G$ with the order induced by a compactification $(G,\theta)$ and identify $\widehat G$ with the subgroup
$\widehat \theta(\widehat G)$ of $\mathbb R.$ For $f \in L^1(G),$ we consider its spectrum $\Omega(f) \subset \mathbb R.$ We define the injective $C^*$-algebra homomorphism $\Theta \, : \, C(G) \rightarrow C_b(\mathbb R)$ by $\Theta(f) := f \circ \theta.$
%
%
\begin{lem}\label{lem12}
$\Theta(T(G)) \subset T(\mathbb R),$ $\Theta(C(G)) \subset U(\mathbb R),$
$m(\Theta(f)) = \int_G f \, d\sigma,$ \newline
$\widehat {\Theta(f)} = \widehat f,$ $\Omega(\Theta(f)) = \Omega(f),$
$\Theta(B(G)) \subset BU(\mathbb R),$ $\Theta(IOHA(G)) \subset IHU(\mathbb R).$
\end{lem}
{\bf Proof}
The first assertion follows since $\chi \in \widehat G \Rightarrow \Theta(\chi) \in \widehat R,$ the second then follows since the Stone-Weierstrass theorem implies that $T(G)$ is dense in $C(G),$ and the third follows
from the theorem of averages (\cite{arnold}, p. 286)
and implies the remaining assertions.
%
%
\begin{theo}\label{thm3}
	If $w \in L_{*}^1(G),$ $w \geq 0,$ $\log w \in L^2(G)$ and $\Omega(w)$ is bounded, then the spectral factor $h = S(w),$ whose existence is ensured by Theorem \ref{thm2}, also has a bounded spectrum and $\Omega(h) \subset [0,\tau]$ where $[-\tau, \tau]$ is the smallest interval containing $\Omega(w).$
\end{theo}
{\bf Proof} Let $[-\tau, \tau]$ where $0 \leq 0 < \infty,$ be the smallest interval containing $\Omega(w).$ Theorem 2 gives $h_n \in IOHA(G)$ such that $\lim_{n \rightarrow \infty} ||h - h_n||_2 = 0$ and $\Omega(|h_n|^2) \subset [-\tau,\tau].$
Lemma \ref{lem12} ensures that $\Theta(h_n) \in IHU(\mathbb R)$ and
$\Omega(|\Theta(h_n)|^2) = \Omega(|h_n|^2) \subset [-\tau,\tau].$
Theorem \ref{thm1} then implies that $\Omega(h_n) \subset [0,\tau]$ hence
$\Omega(h) \subset [0,\tau].$
\begin{cor}\label{cor1}
If $w \in C(G),$ $w \geq 0$ and $\Omega(w)$ is bounded
then $\log w \in L^2(G).$ Therefore the conclusion of Theorem \ref{thm3} holds.
\end{cor}
{\bf Proof} Clearly $w \in C(G) \subset L_{*}^1(G)$ and $f := \Theta(w) \in U(\mathbb R)$ and $\Omega(f)$ is bounded. In (\cite{lawton2}, Corollary 1.1) we proved that $\log f \in B^p(\mathbb R),$ the Besicovitch \cite{besicovitch1}, \cite{besicovitch2} space of almost periodic functions for all $p \geq 1.$ Therefore (\cite{lawton2}, Lemma 2.2) implies that $w \in L^2(G),$ so the hypothesis and conclusion of Theorem \ref{thm3} holds.
\begin{remark}\label{rem4}
Theorem \ref{thm3} extends results in \cite{lawton1} for special archimedian orders on $\widehat {\mathbb T^2} = \mathbb Z^2.$
\end{remark}
\begin{cor}\label{cor2}
If $h \in C(G)$ has a bounded spectrum $\Omega(h) \subset [0,\infty)$
and $F \in HOL(\mathbb U)$ is the extension of $\Theta(h),$
then $h$ is outer iff $F$ has no zeros in $\mathbb U^o.$
\end{cor}
{\bf Proof} We may assume that $h$ is not the zero function. Let $\tau \geq 0$ be the smallest number such that $\Omega(h) \subset [0,\infty).$ Define $w := |h|^2.$ Then $[-\tau,\tau]$ is the smallest closed interval with $\Omega(w) \subset [-\tau, \tau].$
Let $h_o := S(w)$ be the spectral factor of $w$ ensured by Theorem \ref{thm2}. If $h$
is outer then $h = h_o$ so the proof of Theorem \ref{thm3} implies that there exists $h_n \in IOHA(G)$ with $\lim_{n \rightarrow \infty} ||h - h_n||_2 = 0.$
Then $\Theta(h), \Theta(h_n)$ have bounded spectra so extend to entire
functions $H, H_n,$ respecively. Lemma 3 implies that $H_n$ has no zeros in $\mathbb U^o.$ Furthermore
$H_n$ converges to $H$ uniformly on compact subsets of $\mathbb U^o,$ therefore a theorem of Hurwitz \cite{hurwitz} implies that $h$ has no zeros in $\mathbb U^o.$ To prove the converse assume that $H$ has no zeros in $\mathbb U^o.$ Then $e^{-{\tau}{2}z}H(z)$ is Ahiezer's spectral factor of the entire function $H(z)\overline {H(\overline z)}$ which equals the entire function $H_o(z)\overline {H_o(\overline z)}$ since they both equal $\Theta(w)$
on $\mathbb R.$ Theorem \ref{thm3} then implies that $H = H_o$ hence $h = h_o$ so $h$ is outer.
\begin{remark}\label{rem5}
Corollary \ref{cor2} proves our previous Conjecture (\cite{lawton2}, Conjecture 2.1).
\end{remark}
{\bf Acknowledgments} The author thanks Professor August Tsikh for insightful discussions
and colleagues in Austria and Thailand for sending references required to write this paper.

\Finishall

\begin{thebibliography}{10}
\bibitem{ahiezer} N. I. Ahiezer, {\it On the theory of entire
functions of finite degree}, Doklady Akad. Nauk SSSR, 63 (1948)
475--478, MR 10,289.

\bibitem{arens} R. Arens, {\it The group of invertible elements
in a commutative Banach algebra}, Studia Math., 1 (1963) 21--23.

\bibitem{arnold} V. I. Arnold, {\it Mathematical Methods of Classical Mechanics},
Springer, New York, 1978.

\bibitem{besicovitch1} A. S. Besicovitch,\textit{On generalized almost periodic functions},
Procedings of the London Mathematical Society, 25(2) (1926) 495--512.

\bibitem{besicovitch2} A. S. Besicovitch, {\it Almost Periodic
Functions}, Cambridge University Press, 1954.

\bibitem{beurling} A. Beurling, {\it On two problems concerning
linear transformations Hilbert space}, Acta Math. 81 (1949) 239--255.

\bibitem{boas1} R. P. Boas, \textit{\it Functions of exponential
type. I}, Duke Math. J., 11 (1944) 9--15.

\bibitem{boas2} R. P. Boas, {\it Entire Functions}, Academic Press,
New York, 1954.


\bibitem{bohr1} H. Bohr, {\it Zur Theorie der fastperiodischen
Funktionen I}. Acta Mathematica, 45 (1924) 29--127.

\bibitem{bohr2} H. Bohr, {\it Zur Theorie der fastperiodishen Functionen
III Teil. Dirichletenwich-lung analytisher Functionen}, Acta Math. 47 (1926) 237-281.

\bibitem{bohr3} H. Bohr, {\it Almost Periodic Functions}, Chelsea, New York, 1951.

\bibitem{fejer} L. Fej\'{e}r, {\it \"{U}ber trigonometricshe Polynome},
J. Reine und Angewandte Mathematik, 146 (1915) 53--82.

\bibitem{deng} G. Deng, {\it Integral representations of harmonic functions in half spaces},
Bulletin des Sciences Math\'{e}matiques, 131 (2007) 53--59.

\bibitem{gelfand1} I. M. Gelfand, {\it Normierte rings},
Mat. Sbornik, 9 (1941) 3--24.

\bibitem{gelfand2} I. M. Gel'fand et al, {\it Generalized Functions, vols I–-VI},
Academic Press, 1964. (Translated from Russian.)

\bibitem{hausdorff} F. Hausdorff, {\it Eine Ausdehnung des Parsevalschen Satzes
\"{u}ber Fourierreihen}, Mathematische Zeitschrift, 16 (1923) 163–169.

\bibitem{heins} M. H. Heins, {\it Entire functions with bounded minimum modulus;
subharmonic functions analogues}, Ann. Math. 49(2) (1948) 200--213.

\bibitem{helsonlow1} H. Helson and D. Lowdenslager,
\textit{Prediction theory and Fourier series in several variables},
Acta Math., 99(1) (1958) 165--202.

\bibitem{hurwitz} A. Hurwitz, {\it \"{U}ber die Bedingungen, unter welchen eine
Gleichung nur W\"{u}zeln mit reelen Teilenbesitz}, Math. Ann., 46 (1989) 273--284.

\bibitem{jensen} J. Jensen, {\it Sur unnouvelet important th\'{e}or\'eme de la
th\'{e}orie des fonctions}, Acta Mathematica, 22 (1899) 359--364.

\bibitem{laugesen} R. S. Laugesen, {\it Harmonic Analysis Lecture Notes},
University of Illinois at Urbana Champaign, Richard S. Laugeson, 2017.
arxiv:0903.3845v2

\bibitem{lawton1} W. Lawton, {\it Spectral factorization of trigonometric
polynomials and lattice geometry}, Acta Arithmetica, 155(3) (2012) 311--321.
arxiv.org/abs/1208.5590

\bibitem{lawton2} W. Lawton, {\it Distribution of small values of Bohr
almost periodic functions with bounded spectrum}, Journal of the Siberian 
Federal University to appear {\bf 12}(6) December (2019). arxiv.org/abs/1904.09373

\bibitem{levi} F. W. Levi, {\it Ordered groups}, Proceedings
of the Indian Academy of Sciences, A16 (1942) 256--263.

\bibitem{levin} B. Ja. Levin, {\it Distribution of Roots of Entire Functions},
Translations of Mathematical Monographs, Volume 5, Revised Edition, American
Mathematical Society, 1964.

\bibitem{mahler1} K. Mahler, {\it An application of Jensen's formula to polynomials},
Mathematica, 7 (1960) 98--100.

\bibitem{mahler2} K. Mahler, {\it On some inequalities for polynomials
in several variables}, J. London Mathematical Society, 37 (1962) 341--344.

\bibitem{paley} R. E. A. C. Paley and N. Wiener, {\it Fourier Transforms in the
Complex Domain}, Amer. Math. Soc. Coll. Publ., Vol. 19, New York, 1934.

\bibitem{pontryagin} L. Pontryagin, {\it Topological Groups},
Princeton University Press, 1946.

\bibitem{ribenboim} P. Ribenboim, {\it The Theory of Classical Valuations},
Springer Monographs in Mathematics, New York, 1999.

\bibitem{riesz1}  F. Riesz, {\it \"{U}ber ein Problem des Hern Carath\'{e}odory},
J. Reine und Angewandte Mathematik, 146 (1915) 83--87.

\bibitem{riesz2}  F. Riesz and B. Sz.-Nagy, {\it Functional Analysis},
Frederick Ungar Publishing Company, New York, 1955.

\bibitem{royden} H. L. Royden, {\it Function algebras},
Bulletin of the American Mathematical Society, 69 (1963) 281--298.

\bibitem{rudin1} W. Rudin, {\it Fourier Analysis on Groups},
John Wiley \& Sons, New York, 1987.

\bibitem{rudin2} W. Rudin, {\it Real and Complex Analysis},
McGraw-Hill, New York, 1990.

\bibitem{schwartz1} L. Schwartz, {\it Th\'{e}orie des Distributions,1-1}, Hermann, 1951.

\bibitem{schwartz2} L. Schwartz, {\it Transformation de Laplace des distributions},
Comm. S\'{e}m. Math. Univ. Lund, Tome Suppl\'{e}mentaire, (1952) 196--206.

\bibitem{szego1} G. Szeg\"{o}, {\it \"{U}ber die Randwerte
analytischer Funktionen}, Math. Ann. 84 (1921) 232--244.

\bibitem{szego2} G. Szeg\"{o}, {\it Zur theorie der fastperiodischen Funktionen},
Math. Ann. 96 (1926) 378--382.

\bibitem{wiener} N. Wiener, {\it Tauberian theorems}, Ann. of Math., 33 (2) (1932) 1--100.

\bibitem{young} W. H. Young, {\it On the determination of the summability of a
function by means of its Fourier constants}, Proc. London Math. Soc., 12 (1913) 71--88.

\end{thebibliography}
\end{document}